\documentclass[12pt]{amsart}

\usepackage[dvips]{graphicx}
\usepackage{amsfonts,amssymb,amscd,bm}
\newcommand{\R}{{\mathbb R}}

\newcommand{\Z}{{\mathbb Z}}

\newcommand{\T}{{\mathbb T}}

\newcommand{\e}{{\varepsilon}}
\newcommand{\al}{{\alpha}}
\newcommand{\be}{{\beta}}

\newcommand{\De}{{\Delta}}
\newcommand{\si}{{\sigma}}
\newcommand{\Si}{{\Sigma}}

\newcommand{\ti}{\tilde}

\newcommand{\bu}{\bullet}

\newcommand{\id}{{\mathrm{id}}}
\newcommand{\pr}{{\mathrm{pr}}}

\newcommand{\lk}{{\mathrm{lk}}}
\newcommand{\tg}{{\mathrm{TG}}}

\newcommand{\xy}{{\mathrm{XY}}}

\newcommand{\xa}{{\mathrm{X}}}

\newcommand{\za}{{\mathrm{Z}}}
\newcommand{\inp}{{\mathrm{Int}\,}}

\newcommand{\bd}{{\partial}}

\newcommand{\ed}{{\hfill $\blacksquare$}}

\newcommand{\ov}[1]{ \overline{ #1 } }

\newcommand{\SB}{\mathrm{SB}}

\newcommand{\cl}{\mathrm{CL}}

\newcommand{\gd}{\mathrm{GD}}
\newcommand{\lcm}{\mathrm{lcm}}

\newcommand{\zt}{\mathrm{ZT}}
\newcommand{\ra}{\mathrm{RA}}
\newcommand{\tc}{\mathrm{TC}}

\newcommand{\rot}{\mathrm{rot}}

\newcommand{\dxy}{\mathrm{D}_{xy}}

\newcommand{\axz}{\mathrm{A}_{xz}}

\newcommand{\cubicsn}{ \begin{picture}(8,10)(0,2)
  \put(3,0){\line(0,1){10}}
  \put(-1,5){\oval(8,8)[br]}
  \put(7,5){\oval(8,8)[tl]}
 \end{picture} }

\newcommand{\quadrup}{ \begin{picture}(12,10)(0,2)
  \put(3,0){\line(2,5){4}}
  \put(0,3){\line(5,2){10}}
  \put(7,0){\line(-2,5){4}}
  \put(0,7){\line(5,-2){10}}
 \end{picture} }

\newcommand{\tang}{ \begin{picture}(10,10)(0,2)
  \put(1,5){\oval(8,8)[r]}
  \put(9,5){\oval(8,8)[l]}
 \end{picture} }
\newcommand{\tangint}{ \begin{picture}(12,10)(0,2)
  \put(0,3){\line(2,1){10}}
  \put(1,5){\oval(8,8)[r]}
  \put(9,5){\oval(8,8)[l]}
 \end{picture} }
\newcommand{\tangver}{ \begin{picture}(8,8)
  \put(6,4){\oval(8,8)[l]}
  \put(2,4){\circle*{2}}
 \end{picture} }

\newcommand{\trip}{ \begin{picture}(11,10)(0,2)
  \put(0,0){\line(1,1){10}}
  \put(5,0){\line(0,1){10}}
  \put(0,10){\line(1,-1){10}}
 \end{picture} }

\newcommand{\tripver}{ \begin{picture}(11,10)(0,2)
  \put(0,0){\line(1,1){10}}
  \put(0,5){\line(1,0){10}}
  \put(0,10){\line(1,-1){10}}
  \put(5,5){\circle*{3}}
 \end{picture} }

\newtheorem{theorem}{Theorem}[section]
\newtheorem{corollary}[theorem]{Corollary}
\newtheorem{proposition}[theorem]{Proposition}
\newtheorem{lemma}[theorem]{Lemma}

\theoremstyle{definition}
\newtheorem{definition}[theorem]{Definition}
\newtheorem{example}[theorem]{Example}

\title{Recognizing trace graphs of closed braids}

\author[fiedler]{T.~Fiedler}
\address{ Laboratoire Emile Picard, Universit\'e Paul Sabatier,
118 route Narbonne, 31062 Toulouse, France}
\email{ fiedler@picard.ups-tlse.fr }

\author[kurlin]{V.~Kurlin}
\address{ Department of Mathematical Sciences,
Durham University,
Durham DH1 3LE, United Kingdom}
\email{ vitaliy.kurlin@durham.ac.uk }

\date{ August 20, 2008, the last version
 is available on www.durham.ac.uk/$\sim$dma0vk}

\begin{document}


\begin{abstract}
To a closed braid in a solid torus we associate a trace graph
 in a thickened torus in such a way that closed braids
 are isotopic if and only if their trace graphs
 can be related by trihedral and tetraherdal moves.
For closed braids with a fixed number of strands,
 we recognize trace graphs up to isotopy and trihedral moves 
 in polynomial time with respect to the braid length.
\end{abstract}

\maketitle



\section{Introduction}

\subsection{Motivation and summary}
\noindent
\smallskip

There is still no efficient solution to the conjugacy problem 
 for braid groups $B_n$ on $n\geq 5$ strands,
 ie with a polynomial complexity in the braid length.
Very promising steps towards a polynomial solution
 were made by Birman, Gebhardt, Gonz\'alez-Meneses \cite{BGG} and Ko, Lee \cite{KL}.
A clear obstruction is that the number of different conjugacy classes
 of braids grows exponentially even in $B_3$, see Murasugi \cite{Mur}.
\smallskip

The conjugacy problem for braids is equivalent to 
 the isotopy classification of closed braids in a solid torus.
To a closed braid in a solid torus we associate 
 a 1-parameter family of closed braids, which
 can be encoded by a new combinatorial object, 
 the labelled \emph{trace graph} in a thickened torus.
\smallskip

We establish the higher order Reidemeister Theorem for closed braids:
\emph{trace graphs determine families of isotopic closed braids
 if and only if they can be related by a finite sequence of
 the trihedral and tetrahedral moves shown in 
 Figure~\ref{fig:TrihedralMove} and Figure \ref{fig:TetrahedralMoves}}, 
 see Theorem~\ref{thm:TraceGraphsMoves}.
We recognize trace graphs of closed braids up to isotopy in a thickened torus
 and trihedral moves in polynomial time with respect to the braid length, 
 see Theorem~\ref{thm:RecognizingUpToTrihedralMoves}.
This is one of very few known polynomial algorithms 
 recognizing complicated topological objects up to isotopy.
\smallskip


\subsection{Basic definitions of braid theory}
\noindent
\smallskip

\begin{figure}[!h]
\includegraphics[scale=1.0]{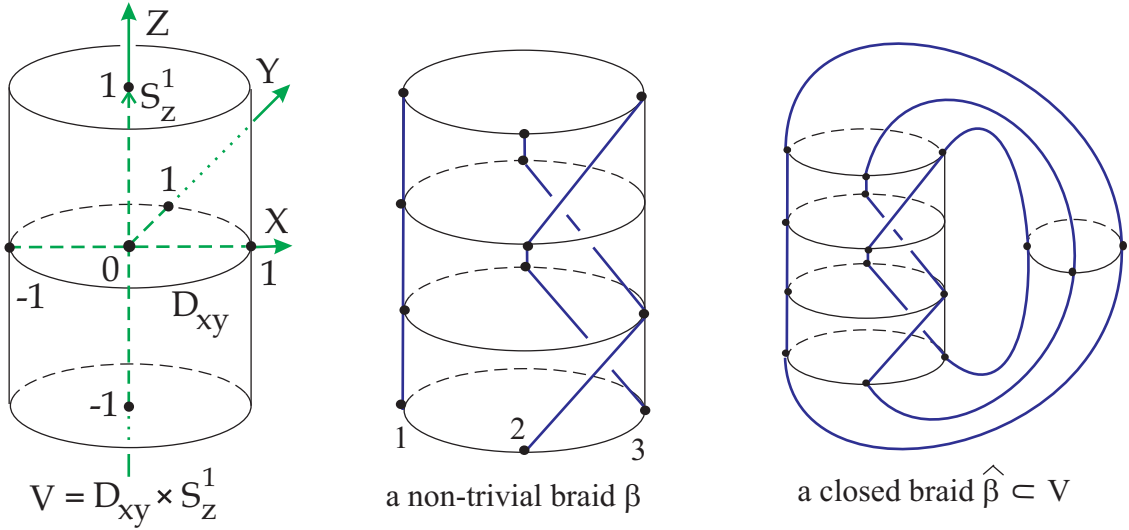}
\caption{A braid and its closure in the solid torus $V$}
\label{fig:BraidClosedBraid}
\end{figure}

We work in the $C^{\infty}$-smooth category.
To explain important constructions we may draw 
 piecewise linear pictures that can be easily smoothed.
Fix Euclidean coordinates $x,y,z$ in $\R^3$.
Denote by $\dxy$ the unit disk at 
 the origin $0$ of the horizontal plane $\xy$.
Introduce the \emph{solid torus}
 $V=\dxy\times S_z^1$, where the oriented circle $S_z^1$
 is the segment $[-1,1]_z$ with the identified endpoints, 
 see the left picture of Figure~\ref{fig:BraidClosedBraid}.
\smallskip

\begin{definition}
\label{def:ClosedBraids}
Mark $n$ points $p_1,\dots,p_n\in\dxy$.
A \emph{braid} $\be$ on $n$ \emph{strands} is
 the image of a smooth embedding of
 $n$ segments into $\dxy\times[-1,1]_z$ such that
\smallskip

\noindent
$\bu$
 the strands of $\be$ are monotonic
 with respect to $\pr_z:\be\to S_z^1$ (see Figure~\ref{fig:BraidClosedBraid});
\smallskip

\noindent
$\bu$
 the lower and upper endpoints of $\be$ are
 $\cup(p_i\times\{-1\})$, $\cup(p_i\times\{1\})$, respectively.
\medskip

\noindent
Identifying the bases $\dxy\times\{z=\pm 1\}$,
 the cylinder $\dxy\times[-1,1]_z$ is converted into the solid torus $V=\dxy\times S_z^1$,
 while a braid $\be\subset\dxy\times[-1,1]_z$ becomes 
 the \emph{closed} braid $\hat\be\subset V$, 
 see the right picture of Figure~\ref{fig:BraidClosedBraid}.
\ed
\end{definition}

\begin{definition}
\label{def:BraidGroup}
Braids are considered up to an \emph{isotopy}, a smooth deformation of 
 the cylinder $\dxy\times[-1,1]_z$, fixed on its boundary.
The equivalence classes of braids
 form the group denoted by $B_n$.
The \emph{product} of braids $\be_1,\be_2$ is the braid $\be_1\be_2$ obtained
 by attaching a cylinder containing $\be_2$ over a cylinder containing $\be_1$.
The \emph{trivial} braid consists of $n$ vertical
 straight segments $\sqcup_{i=1}^n(p_i\times[-1,1]_z)$.
\ed
\end{definition}

The braid group $B_n$ is generated by elementary braids $\si_i$,
 $i=1,\dots,n-1$, where $\si_i$ is a right half-twist of strands $i,i+1$,
 the remaining strands are vertical.
The braid $\be$ in the middle picture of Figure~\ref{fig:BraidClosedBraid}
 can be considered as a product of 3 braids, the lower one is $\si_2$,
 the middle one is trivial, the upper one is $\si_2$, ie $\be=\si_2^2$.
Any braid induces a \emph{permutation} of its endpoints, 
 eg the braid in the middle picture of Figure~\ref{fig:BraidClosedBraid}
 induces the trivial permutation on endpoints 1, 2, 3.
Such a braid is called \emph{pure}.
The closure of any pure braid $\be\in B_n$ consists of $n$ components.


\subsection{Trace graphs of closed braids}
\noindent
\smallskip

\begin{figure}[!h]
\includegraphics{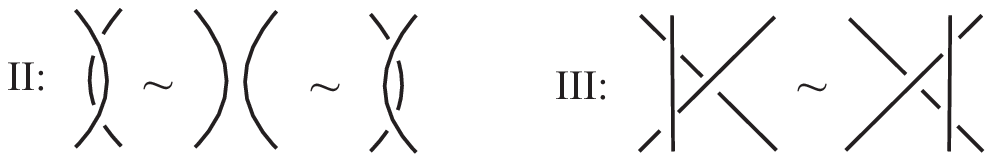}
\caption{Reidemeister moves on braids}
\label{fig:ReidemeisterMovesBraids}
\end{figure}

Closed braids are usually represented by plane diagrams with double crossings.
A classical approach to the isotopy classification of closed braids
 is to use diagram invariants, ie functions defined on 
 plane diagrams and invariant under Reidemester moves II, III in 
 Figure~\ref{fig:ReidemeisterMovesBraids}.
A 1-parameter approach proposed by Fiedler and Kurlin \cite{FK2}
 is to consider the 1-parameter family of diagrams of 
 braids rotated around the core of the solid torus $V$.
This family contains more combinatorial information about a closed braid
 than just one plane diagram and involves such features of braids
 as meridional \emph{trisecants}, straight lines meeting a braid in 3 points
 and contained in a meridional disk $\dxy\times\{z\}$ of the solid torus $V$.

A \emph{long knot} in $\R^3$, a single curve approaching the vertical axis $\za$ at $\pm\infty$,
 can be also rotated in $\R^3$ around $\za$, but closed braids are more naturally rotated in $V$.
It is essential to work in the solid torus instead of $\R^3$ since 
 our 1-parameter family respresents a non-trivial rational homology class 
 in the space of all diagrams. 
A.~Hatcher has proven that the space of diagrams of a prime knot in $\R^3$
 has a finite fundamental group \cite{H}. 
Consequently, its rational first homology group vanishes.

\begin{figure}[!h]
\includegraphics[scale=1.0]{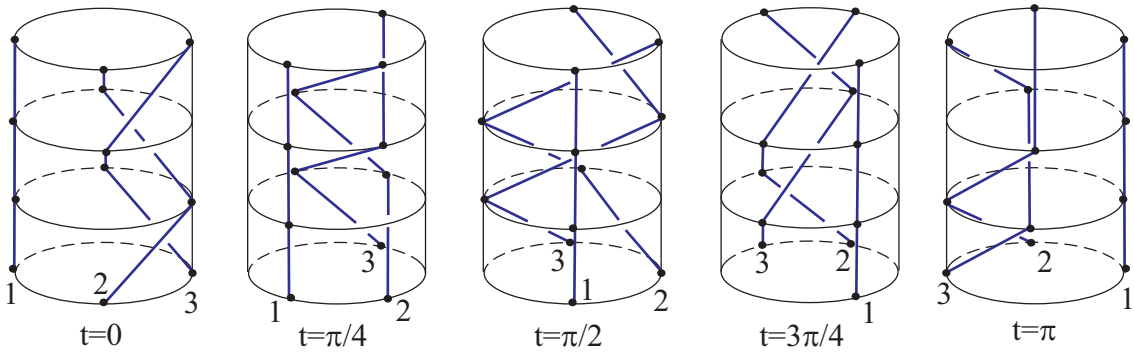}
\caption{Diagrams of rotated braids $\rot_t(\be)$ for the braid $\be$ in Figure~\ref{fig:BraidClosedBraid}}
\label{fig:RotatedBraids}
\end{figure}

\begin{definition}
\label{def:TraceGraphClosedBraid}
Given a closed braid $\hat\be\subset V$ in a general position
 (see more details in subsection~\ref{subs:SingularitiesBraids}),
 consider \emph{rotated} braids $\rot_t(\hat\be)\subset V$ obtained 
 by the rotation of $\hat\be$ through an angle $t\in[0,2\pi)$.
Project each of the rotated braids $\rot_t(\hat\be)$ to 
 the fixed annulus $A_{xz}=[-1,1]_x\times S_z^1\subset V$,
 see Figure~\ref{fig:RotatedBraids}.
The crossings of the resulting diagrams form 
 the \emph{trace graph} $\tg(\hat\be)$ that lives
 in the \emph{thickened} torus $\T=A_{xz}\times S_t^1$, where 
 the time circle $S_t^1$ is $[0,2\pi]$ with the identified endpoints. 
\ed
\end{definition}

\begin{figure}[!h]
\includegraphics[scale=1.0]{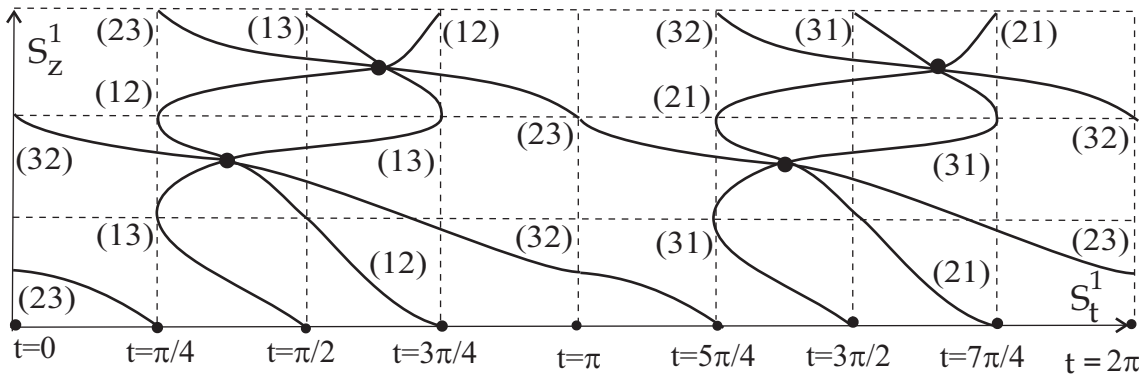}
\caption{The trace graph of the closed braid $\hat\be$ in Figure~\ref{fig:BraidClosedBraid}}
\label{fig:TraceGraphBraid}
\end{figure}

Label arcs of a pure braid $\be\in B_n$ by $1,2,\dots,n$ as 
 in the middle picture of Figure~\ref{fig:BraidClosedBraid}.
Any general point $p$ of the trace graph $\tg(\hat\be)\subset\T$ is 
 a crossing of arcs $i,j$ in the diagram of a rotated braid $\rot_t(\hat\be)$,
 ie the point $p$ evolves in $\T$ following a trace of crossings in the diagrams.
Label the point $p$ by the ordered pair $(ij)$ if the arc $i$ is over the arc $j$
 in the diagram of $\rot_t(\hat\be)$ and by the ordered pair $(ji)$ otherwise.
For non-pure braids, other well-defined markings will be introduced in 
 Definition~\ref{def:TraceCircles}.
The trace graph maps to itself under the time shift $t\mapsto t+\pi$,
 each label $(ij)$ reverses to $(ji)$.
Notice that each labelled closed loop of $\tg(\hat\be)$ is monotonic with respect to 
the vertical circle $S_z^1$, but not with respect to the time circle $S_t^1$.

The trace graph of the piecewise linear closed braid $\hat\be$ in 
 Figure~\ref{fig:BraidClosedBraid} is projected to the torus $\zt=S_z^1\times S_t^1$ 
 and is shown in Figure~\ref{fig:TraceGraphBraid}.
A vertical section $\tg(\hat\be)\cap(\axz\times\{t\})$ of a trace graph consists
 of finitely many points, which are crossings of the diagram of $\rot_t(\hat\be)$,
 eg the zero section $\tg(\hat\be)\cap(\axz\times\{0\})$ contains 2 points
 associated to the crossings of the original braid $\hat\be$.
The section $\tg(\hat\be)\cap(\axz\times\{\pi/4\})$ has 
 2 \emph{tangent} vertices, when the rotated braid $\rot_{\pi/4}(\hat\be)$  
 has 2 simple tangencies (arc 1 over arcs 2 and 3), 
 ie the diagram of $\rot_{\pi/4}(\hat\be)$ changes under Reidemeister moves II.
\smallskip

The braid $\be$ has two meridional trisecants associated to 
 two \emph{triple} vertices of $\tg(\hat\be)$.
Under the rotation of $\be$ through some $t\in(\pi/4,\pi/2)$ and $t\in(\pi/2,3\pi/4)$, 
 the trisecants become perpendicular to the plane of projection, so
 triple intersections appear in the corresponding diagrams of $\rot_t(\hat\be)$.
Around these singular moments the diagrams change under Reidemeister moves III,
 notice that the labels don't change at triple points, see more details about 
 singularities and general position of braids in 
 subsection~\ref{subs:SingularitiesBraids}.
A given closed braid can be reconstructed from its trace graph with labels, see also 
 combinatorial constructions of a trace graph in subsection~\ref{subs:Constructions}.

\begin{figure}[!h]
\includegraphics[scale=1.0]{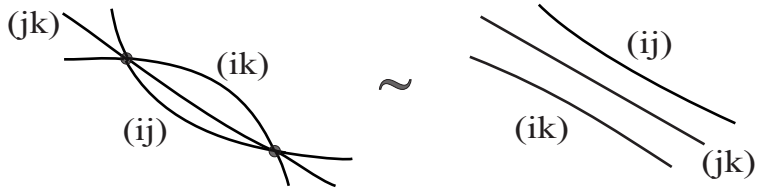}
\caption{Trihedral move on trace graphs}
\label{fig:TrihedralMove}
\end{figure}

\begin{figure}[!h]
\includegraphics[scale=1.0]{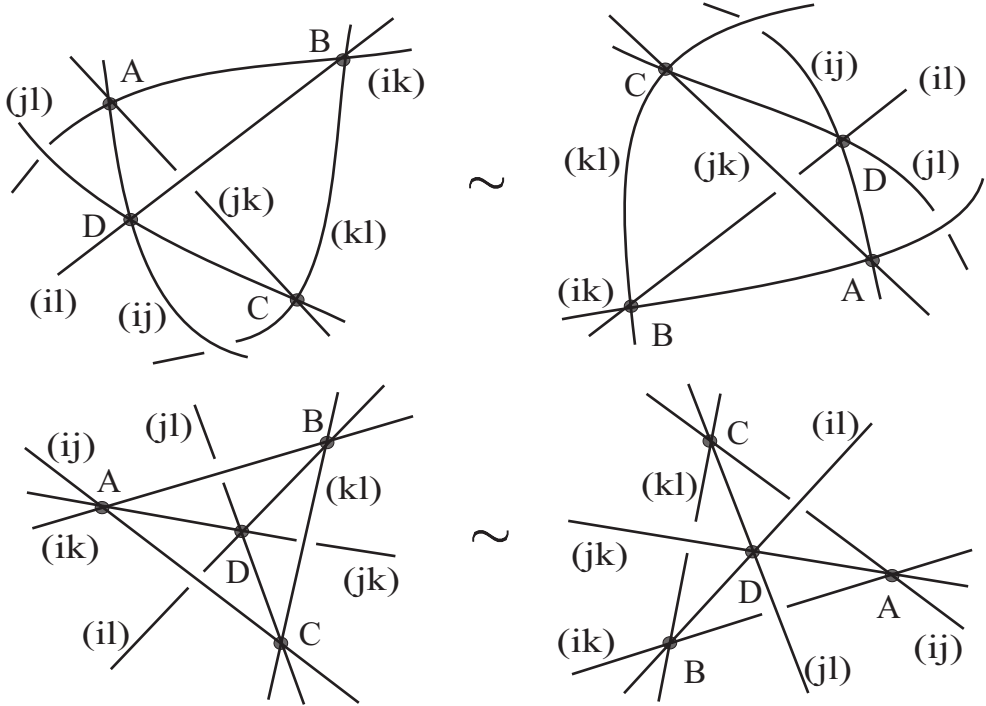}
\caption{Tetrahedral moves on trace graphs}
\label{fig:TetrahedralMoves}
\end{figure}

\begin{theorem}
\label{thm:TraceGraphsMoves}
Closed braids $\hat\be_0,\hat\be_1$ are isotopic in the solid torus $V$
 if and only if their labelled trace graphs $\tg(\hat\be_0),\tg(\hat\be_1)\subset\T$
 can be obtained from each other by an isotopy in $\T$
 and a finite sequence of moves in Figure~\ref{fig:TrihedralMove} 
 and Figure~\ref{fig:TetrahedralMoves}.
\end{theorem}

The trihedral move is associated to the singular situation
 in the space of all closed braids, when the path of rotated braids
 touches the singular subspace of triple intersections $\trip$, ie 
 under the rotation 3 crossings approach each other as in Reidemeister move III, 
 but then go back in the reverse direction without completing Reidemeister move III.
The tetrahedral move is associated to passing through the singular subspace
 of quadruple intersections $\quadrup$, when a 1-skeleton of some tetrahedron
 collapses in the trace graph to a point and then blows up again in a symmetric form.
\smallskip

Theorem~\ref{thm:TraceGraphsMoves} can be used to construct invariants 
 of closed braids reflecting such geometric features as meridional trisecants.
Similar easily computable lower bounds on the number of fiber quadrisecants 
 in knot isotopies were found by Fiedler and Kurlin \cite{FK1}.
On the other hand trace graphs turned out to be complicated topological objects
 that can be recognized up to isotopy in a polynomial time.

\begin{theorem}
\label{thm:RecognizingUpToTrihedralMoves}
Let $\be,\be'\in B_n$ be braids of length $\leq l$.
There is an algorithm of complexity $C(n/2)^{n^2/8}(6l)^{n^2-n+1}$
 to decide whether $\tg(\hat\be)$ and $\tg(\hat\be')$
 are related by isotopy in $\T$ and trihedral moves, 
 the constant $C$ does not depend on $l$ and $n$.
In the case of pure braids, the power $n^2/8$ can be replaced by 1.
If the closure of a braid is a \emph{knot}, a single circle in the solid torus,
 then the complexity reduces to $Cn(6l)^{n-1}$.
\end{theorem}

\noindent
{\bf Acknowledgements.}
The second author is especially grateful to Hugh Morton
 for fruitful suggestions.
He also thanks 
 M.~Kazaryan, V.~Vassiliev 
 for useful discussions.


\section{Studying closed braids in terms of their trace graphs}
\label{sect:ClosedBraidsTraceGraphs}


\subsection{Singularities and general position of closed braids}
\label{subs:SingularitiesBraids}
\noindent
\smallskip

Here we give an outline of the proof of Theorem~\ref{thm:TraceGraphsMoves},
 which immediately follows from a more general result by 
 Fiedler and Kurlin \cite[Theorem~1.4]{FK2} on links in the solid torus $V$
 that can have extrema of the projection to the core $S_z^1$ of $V$.
\smallskip

Codimension~1 singularities of closed braids with respect to the plane projection
 are tangencies of order~1 $\tang$ and triple intersections $\trip$ associated to 
 Reidemeister moves II and III, respectively, see Figure~\ref{fig:ReidemeisterMovesBraids}.   
The Reidemeister theorem says that any isotopy in the space $\SB$ of all closed braids
 (with respect to the Whitney topology) can be approximated by a path 
 transversal to the singular subspace $\Si_{\tang}\cup\Si_{\trip}\subset\SB$.
We extend this approach to 1-parameter families of rotated closed braids.
\smallskip

Codimension~2 singularities of plane diagrams of closed braids  
 are quadruple points $\quadrup$, tangent triple points $\tangint$ and
 tangencies of order~2 $\cubicsn$.
A closed braid  $\hat\be\subset V$ can be put in a \emph{general position} 
 such that the \emph{canonical} loop of 
 rotated braids $\{\rot_t(\hat\be)\}\subset\SB$ is 
 transversal to the codimension~1 subspace $\Si_{\tang}\cup\Si_{\trip}\subset\SB$ and 
 avoids the codimension~2 subspace 
 $\Si_{\quadrup}\cup\Si_{\tangint}\cup\Si_{\cubicsn}\subset\SB$.
Similarly any isotopy of closed braids can be approximated by a path $\{\hat\be_s\}_{s=0}^{s=1}$ 
 such that the cylinder of canonical loops $\{\rot_t(\hat\be_s\}$ is transversal to 
 $\Si_{\quadrup}\cup\Si_{\tangint}\cup\Si_{\cubicsn}\subset\SB$.
Passing through these singularities leads to tetrahedral moves in Figure~\ref{fig:TetrahedralMoves}, 
 trihedral move in Figure~\ref{fig:TrihedralMove} and a move where a triple vertex $\tripver$ 
 of a trace graph passes through a tangent vertex $\tangver$, which doesn't change
 the combinatorial structure of the trace graph with labels, 
 see more details in Fiedler and Kurlin \cite{FK2}.
\smallskip

\begin{figure}[!h]
\includegraphics[scale=1.0]{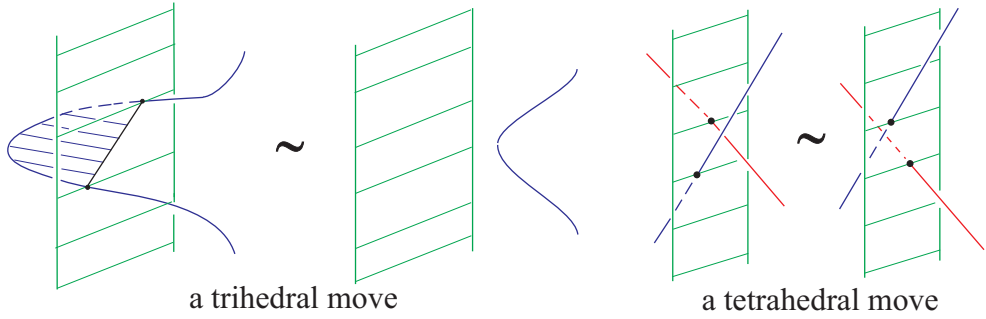}
\caption{A trihedral move and a tetrahedral move for braids.}
\label{fig:MovesOnBraids}
\end{figure}

A geometric interpretation of a trihedral move and tetrahedral move 
 at the level of closed braids is shown in Figure~\ref{fig:MovesOnBraids}.
In a tetrahedral move two arcs intersect a wide band bounded by 
 other two arcs, so two intersection points swap their heights.
The first move in Figure~\ref{fig:TetrahedralMoves} applies when
 the intermediate oriented arcs go together
 from one side of the band to another like $\rightrightarrows$.
The second move in Figure~\ref{fig:TetrahedralMoves} means that the arcs are antiparallel
 as in the British rail mark $\rightleftarrows$.


\subsection{Combinatorial constructions of a trace graph}
\label{subs:Constructions}
\noindent
\smallskip

First we show how to create the trace graph using an algebraic form of a braid.

\begin{lemma}
\label{lem:AlgConstruction}
Let $\be\in B_n$ be a braid of length $l$.
Then the closure $\hat\be$ is isotopic in the solid torus $V$ to a closed braid 
 whose trace graph contains $2l(n-2)$ triple vertices.
\end{lemma}
\begin{proof}
Let $\Delta \in B_n$ be Garside's element \cite{Gar}, i.e.
 $\Delta^2$ is a generator of the centre of $B_n$,
 the full twist of $n$ strands.
The rotation of a braid $\be\in B_n$ can be considered as
 a commutation of $\be$ with $\De^2$.
So the canonical loop of rotated closed braids $\rot_t(\hat\be)$ 
 is represented by the sequence of the closures of the following braids:
$$\beta\to \Delta\Delta^{-1}\beta\to
  \Delta^{-1}\beta\Delta\to
  \Delta^{-1}\Delta\be'\to \be'\to
  \Delta\Delta^{-1}\beta'\to
  \Delta^{-1}\be'\Delta\to
  \Delta^{-1}\De\be\to \be.$$

The first arrow in the sequence consists of Reidemester moves II
 creating couples of symmetric crossings.
The second arrow represents an isotopy of the diagram
 when we push $\De$ through the trivial part of the closed braid $\hat\be$,
 ie we cyclically shift the letters of $\De\De^{-1}\be$ to get $\De^{-1}\be\De$.
The third arrow shows how $\De$ acts on $\be$ from the right.
After we get a new braid $\be'$, we apply
 the same transformation and finish with $\be$
 since $\be\De=\De\be'$ implies that $\be'\De=\De\be$.
\smallskip

For $n=3$, we have $\Delta=\sigma_1\sigma_2\sigma_1$.
We need to consider only the two generators $\si_1,\si_2$
 and their inverses.
We apply braid relations corresponding to
 Reidemeister moves II and III associated to
 tangent and triple vertices of $\tg(\be)$.
$$\begin{array}{l}
\sigma_1\Delta = \sigma_1(\sigma_1\sigma_2\sigma_1) \to
 \sigma_1(\sigma_2\sigma_1\sigma_2) =\Delta\sigma_2,\\
\sigma_2\Delta = \sigma_2(\sigma_1\sigma_2\sigma_1) \to
 (\sigma_1\sigma_2\sigma_1)\sigma_1 = \Delta\sigma_1,\\
\sigma_1^{-1}\Delta = \sigma_1^{-1}(\sigma_1\sigma_2\sigma_1) \to
 \sigma_2\sigma_1\to (\sigma_1\sigma_1^{-1})\sigma_2\sigma_1 \to
 \sigma_1(\sigma_2\sigma_1\sigma_2^{-1}) = \Delta\sigma_2^{-1},\\
\sigma_2^{-1}\Delta = \sigma_2^{-1}(\sigma_1\sigma_2\sigma_1) \to
 (\sigma_1\sigma_2\sigma_1^{-1})\sigma_1 \to \sigma_1\sigma_2 \to
 \sigma_1\sigma_2(\sigma_1\sigma_1^{-1}) = \Delta\sigma_1^{-1}.
\end{array}$$

Notice that the sequence is canonical in the case of
 a generator and almost canonical in the case of
 an inverse generator.
Indeed we can replace the above sequence by
 $\sigma_1^{-1}\Delta \to \Delta\sigma_2^{-1}$ by
 $\sigma_1^{-1}(\sigma_1\sigma_2\sigma_1)\to \sigma_2\sigma_1\to
  \sigma_2\sigma_1(\sigma_2\sigma_2^{-1})\to
  (\sigma_1\sigma_2\sigma_1)\sigma_2^{-1}$.
Pushing $\Delta$ through a generator or its inverse
 creates exactly $n-2$ triple points.
So we end up with $2l(n-2)$ triple points,
 because we push $\Delta$ twice through the braid.
\end{proof}

Now we construct a trace graph of a closed braid in a geometric way.
\smallskip

\begin{lemma}
\label{lem:GeoConstruction}
Let $\be\in B_n$ be a braid of length $l$.
Then the closure $\hat\be$ is isotopic in the solid torus $V$ to a closed braid
 whose trace graph consists of elementary blocks associated
 to the generators (and their inverses) of $B_n$ similar to 
 Figure~\ref{fig:ElementaryTraceGraphs}.
\end{lemma}
\begin{proof}
Figure~\ref{fig:ElementaryTraceGraphs} shows the trace graphs of
 the elements $\si_1$ and $\si_2^{-1}$ in the braid group $B_4$.
In general we mark out the points $\psi_k=2^{1-k}\pi$,
 $k=0,\dots,n-1$ on the boundary of the bases $\dxy\times\{\pm 1\}$.
The $0$-th point $\psi_0=2\pi$ is the $n$-th point.
\smallskip

The crucial feature of the distribution $\{\psi_k\}$
 is that all straight lines passing through two points
 $\psi_j,\psi_k$ are not parallel to each other.
Firstly we draw all strands in the cylinder
 $\bd\dxy\times[-1,1]_z$.
Secondly we approximate with the first derivative
 the strands forming a crossing by smooth arcs,
 see the left pictures of Figure~\ref{fig:ElementaryTraceGraphs}.
\smallskip

\begin{figure}[!h]
\includegraphics[scale=1.0]{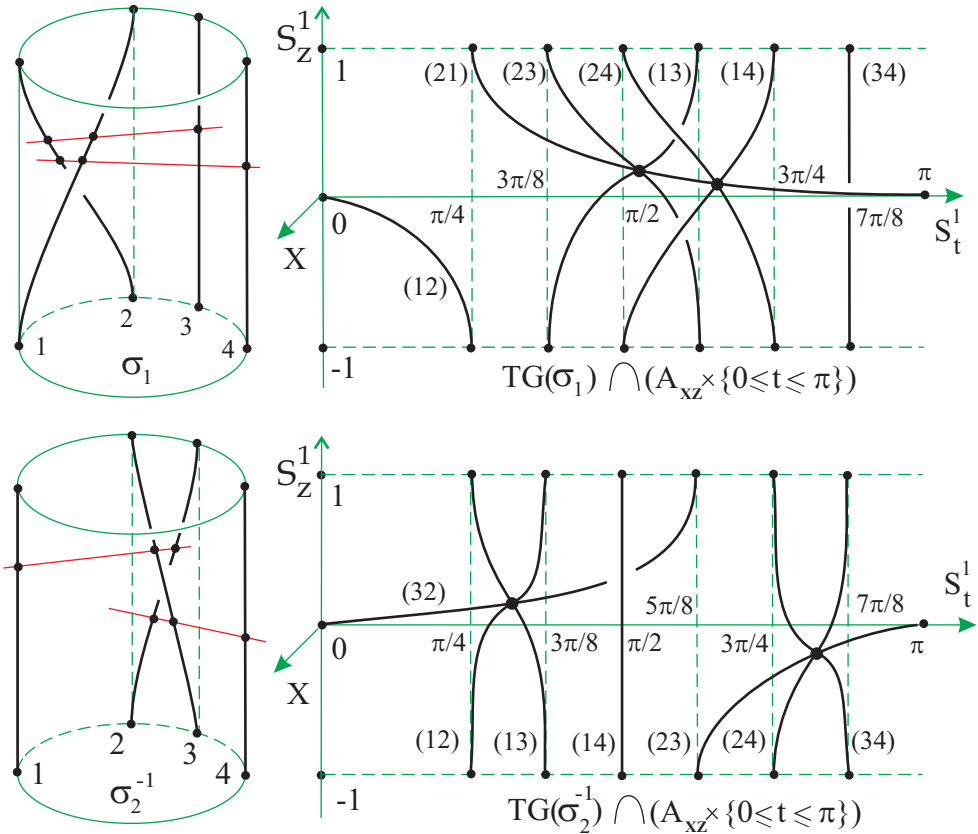}
\caption{Half trace graphs of the 4-braids
$\si_1,\si_2^{-1}\in B_4$.}
\label{fig:ElementaryTraceGraphs}
\end{figure}

Then each elementary braid $\si_i$ constructed as above 
 has exactly $n-2$ meridional trisecants,
 one trisecant through the strands $i,i+1$ and $j$ for each $j\neq i,i+1$.
Each trisecant is associated to a triple vertex of the trace graph,
 see 4 meridional trisecants in the left pictures of Figure~\ref{fig:ElementaryTraceGraphs}.
The right pictures in Figure~\ref{fig:ElementaryTraceGraphs} contain
 the trace graphs of the corresponding 4-braids.
The braids are not in general position, eg parallel strands 3 and 4 
 lead to the vertical arc labelled with $(34)$, but
 we may slightly deform such a braid, 
 which makes the projection $\tg\to S_t^1$ generic.
\end{proof}


\section{Two splittings of trace graphs of closed braids}
\label{sect:TwoSplittings}


\subsection{A trace graph splits into trace circles}
\label{subs:SplittingTraceCircles}
\noindent
\smallskip


\begin{definition}
\label{def:TraceCircles}
A \emph{trace} circle in the trace graph of a closed braid is 
 a closed loop that does not change its direction at triple vertices.
\end{definition}

\begin{lemma}
\label{lem:SplittingTraceCircles}
For a braid $\be\in B_n$, let $(n_1,\dots,n_{m})$
 be the lengths of cycles in the induced
 permutation $\ti\be\in S_n$.
Enumerate all components of the closure
 $\hat\be$ by $1,\dots,m$.
Set $N(\be)=\sum\limits_{i=1}^{m}(n_i-1)+
 2\sum\limits_{i<j}\gcd(n_i,n_j)$,
 gcd is the greatest common divisor.
\smallskip

\noindent
{\bf (a)}
The trace graph $\tg(\hat\be)\subset\T$ splits into
 $N(\be)$ trace circles that can be marked by $T_{(ij)[k]}$, 
 where indices $i,j\in\{1,\dots,m\}$, $k=1,\dots,\gcd(n_i,n_j)$.
\medskip

\noindent
{\bf (b)}
If the braid $\be$ is pure,
 ie the permutation $\ti\be$ is trivial,
 then $N(\be)=n(n-1)$.\\
If the closure $\hat\be$ of the braid $\be$ is a knot,
 then $N(\be)=n-1$.\\
For any braid $\be\in B_n$, we have
 $n-1\leq N(\be)\leq n(n-1)$.
\end{lemma}
\begin{proof}
{\bf (a)}
Denote by $p_1,\dots,p_n$ the intersections of $\hat\be$
 with $\dxy\times\{-1\}$, ordered by the orientation of $\hat\be$.
Suppose that $p_r,p_s$ belong to the $q$-th component of $\hat\be$.
This component corresponds to a cycle of the length $n_q$
 of the permutation $\ti\be\in S_n$.
\smallskip

If we push $p_r,p_s$ along their strands in $\hat\be$,
 the associated point in $\tg(\hat\be)$ goes along
 a trace circle and comes to the point corresponding to
 the next pair (say) $(p_{r+1},p_{s+1})$.
This process continues until we come to the original pair
 $(p_r,p_s)$ after $n_q$ steps along the cycle of $\ti\be$
 having passed through $n_q$ of $n_q(n_q-1)$ ordered pairs.
For each cycle of length $n_q$ of $\ti\be\in S_n$,
 we get $n_q-1$ trace circles.
The trace circles are distinguished by 
 non-zero differences $r-s\pmod{n_q}\in\{1,\dots,n_q-1\}$.
\smallskip

Assume that $p_r,p_s$ are in different components $i\neq j$ 
 of $\hat\be$, associated to cycles of lengths $n_i,n_j$.
Then the process above terminates after
 $\lcm(n_i,n_j)$ steps, $\lcm$ is the lowest common multiple,
 since at each step indices $r,s$ shift by 1 in two sets of lengths $n_i,n_j$.
For any two cycles of lengths $n_i,n_j$ in $\ti\be$,
 we get $2\gcd(n_i,n_j)$ trace circles split
 into pairs symmetric with respect to $t\mapsto t+\pi$.
\medskip

\noindent
{\bf (b)}
If $\hat\be$ is a knot then the permutation
 $\ti\be$ is cyclic, ie $m=1$, $N(\be)=n-1$.
For a pure braid $\be\in B_n$, we have $n_1=\cdots=n_m=1$,
 hence $N(\be)=n(n-1)$.
The upper estimate $N(\be)\leq n(n-1)$ geometrically follows
 from the fact that all trace circles
 are monotonic in the direction $S_z^1$ and
 each meridional disk $\dxy\times\{z\}$ intersects $\hat\be$
 in exactly $n$ points leading to $n(n-1)$ crossings appearing under the rotation.
\smallskip

Let $n_1$ be minimal among all lengths $n_i>1$.
Under the map
 $(n_1,n_2,\dots,n_{m})\mapsto
 (\overbrace{1,\dots,1}^{n_1},n_2,\dots,n_{m})$,
 the number $N(\be)$ of trace circles increases by
$$n_1(n_1-1)+2n_1(m-1)-(n_1-1)-2\sum\limits_{i=2}^m\gcd(n_1,n_i)\geq$$
$$\geq (n_1-1)^2+2n_1(m-1)-2\sum\limits_{i=2}^m n_1=(n_1-1)^2\geq 0.$$
So $N(\be)$ is minimal
 if $\hat\be$ is a knot and maximal if $\be$ is pure.
\end{proof}
\smallskip

\noindent
The trace circles $T_{(ii)[k]}$ are associated to the points of
 $(\dxy\times\{z\})\cap\hat\be$ from
 the $i$-th component of $\hat\be$.
If we index these points by $1,\dots,n_i$
 according to the orientation of the $i$-th component then
 the number $k\in\{1,\dots,n_i-1\}$ in the notation $(ii)[k]$
 is well-defined as the difference between the indices modulo $n_i$.
\smallskip

\noindent
The trace circles $T_{(ij)[k]}$ with $i\neq j$ are generated
 by the $i$-th and $j$-th components of $\hat\be$.
Then $k\in\{1,\dots,\gcd(n_i,n_j)\}$ is defined up to
 cyclic permutation, ie if a trace circle is marked by 1,
 this defines markings on other circles $T_{(ij)[k]},T_{(ji)[k]}$, $k>1$.
\smallskip

\begin{figure}[!h]
\includegraphics[scale=1.0]{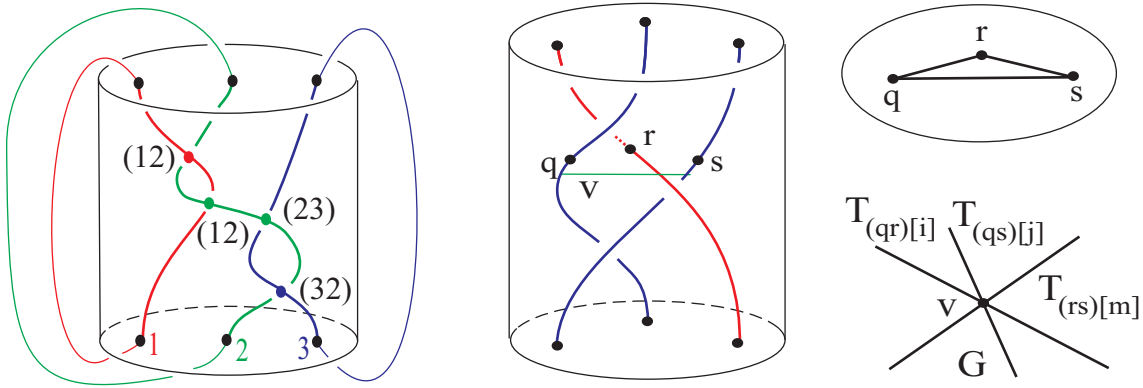}
\caption{Markings of crossings and trace circles.}
\label{fig:MarkingsTraceCircles}
\end{figure}

If $\hat\be$ is pure then the markings reduce to
 ordered pairs $(ij)$, see Figure~\ref{fig:MarkingsTraceCircles}.
If the closure $\hat\be$ is a knot, the markings
 are well-defined numbers $[k]\in\{1,\dots,n-1\}$.

\begin{lemma}
\label{lem:PositionTraceCircles}
Let $\be\in B_n$ be a braid.
Consider trace circles $T_{(qr)[i]},T_{(qs)[j]},T_{(rs)[m]}$
 passing through a triple vertex $v\in\tg(\hat\be)$.
\smallskip

\noindent
{\bf (a)}
The trace circle $T_{(qs)[j]}$ passes between $T_{(qr)[i]}$
 and $T_{(rs)[m]}$ at the vertex v.
\smallskip

\noindent
{\bf (b)}
If $\be$ is pure then each trace circle $T_{(ij)}$ maps to $T_{(ji)}$
 under $t\mapsto t+\pi$.\\
If $\hat\be$ is a knot then each trace circle $T_{[m]}$ maps to
 $T_{[n-m]}$ under $t\mapsto t+\pi$.
\end{lemma}
\begin{proof}
{\bf (a)}
Let the components $\hat\be$ indexed by $q,r,s$ form
 a triple intersection $(qrs)$ associated to
 the vertex $v\in\tg(\hat\be)$, see Figure~\ref{fig:MarkingsTraceCircles}.
Consider a disk $\dxy\times\{z\}$
 slightly above the triple intersection $(qrs)$.
In the disk we see 3 points of the arcs $q,r,s$.
These points form a triangle, the angle at the point of
 $r$ is close to $\pi$.
\smallskip

Denote by $t_{qr},t_{rs},t_{qs}\in S_t^1$ the time moments
 when the corresponding arcs in the closed braid $\hat\be$
 form a crossing under $\pr_{xz}$.
Then $t_{qs}$ is between $t_{qr}$ and $t_{rs}$.
So the crossing $(qs)$ is associated to the middle circle
 $T_{(qs)[j]}$ between $T_{(qr)[i]}$ and $T_{(rs)[m]}$.
\medskip

\noindent
{\bf (b)}
For a pure braid $\be$, let a point $p\in T_{(ij)}$ correspond
 to a crossing $(p_i,p_j)\subset\hat\be\cap(\dxy\times\{z\})$.
Under $t\mapsto t+\pi$, the crossing $(p_i,p_j)$
 converts to the reversed crossing $(p_j,p_i)$ associated to
 the trace circle marked by $(ji)$, see Figure~\ref{fig:SplittingIntoTraceCircles}.
\smallskip

\begin{figure}
\includegraphics[scale=1.0]{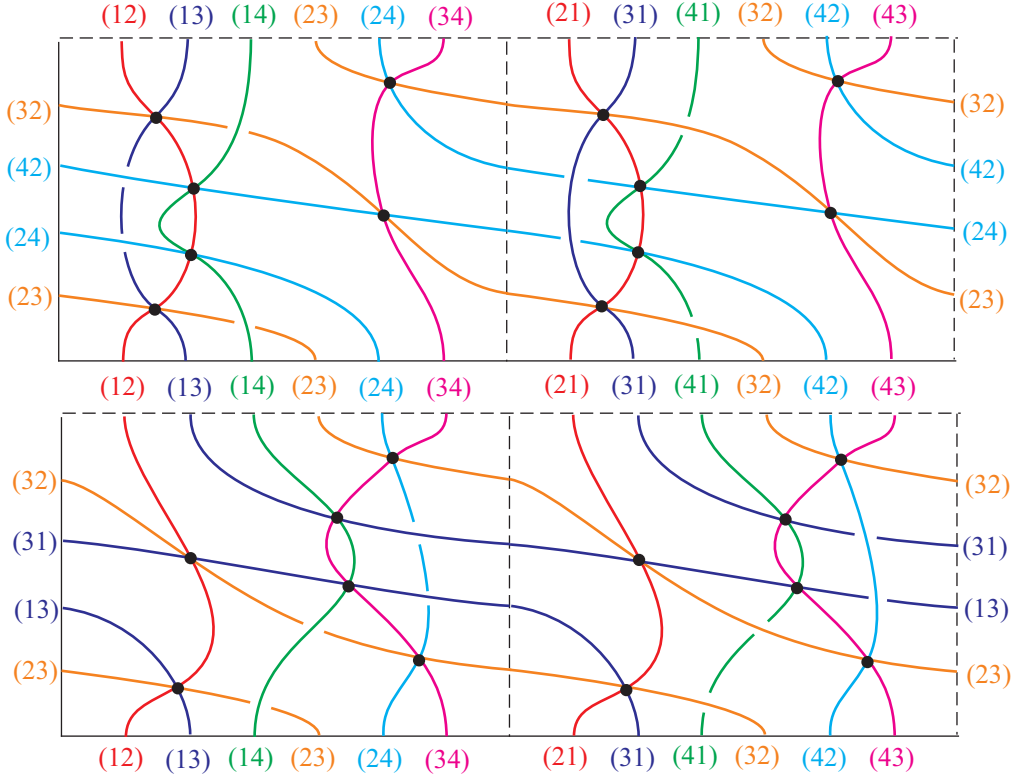}
\caption{Splittings of $\tg(\widehat{\si_2\si_3^2\si_2})$,
 $\tg(\widehat{\si_2\si_1^2\si_2})$ into trace circles .}
\label{fig:SplittingIntoTraceCircles}
\end{figure}

If $\hat\be$ is a knot, order all intersections
 $(p_1,\dots,p_n)\subset\hat\be\cap(\dxy\times\{z\})$
 according to the orientation of $\hat\be$.
Then the marking $[m]$ of a crossing $(p_r,p_s)$ is
 $r-s\pmod{n}$, hence the reversed crossing
 $(p_s,p_r)$ has the marking $s-r\equiv n-m\pmod{n}$.
\end{proof}


\subsection{A trace graph splits into level subgraphs}
\label{subs:SplittingLevelSubgraphs}
\noindent
\smallskip

Here we split the trace graph $\tg(\hat\be)$ of a closed braid $\hat\be$ into
 level subgraphs, trivalent graphs $S^{(k)}$, where $\be\in B_n$ and $k=1,\dots,n-1$.

\begin{figure}[!h]
\includegraphics[scale=1.0]{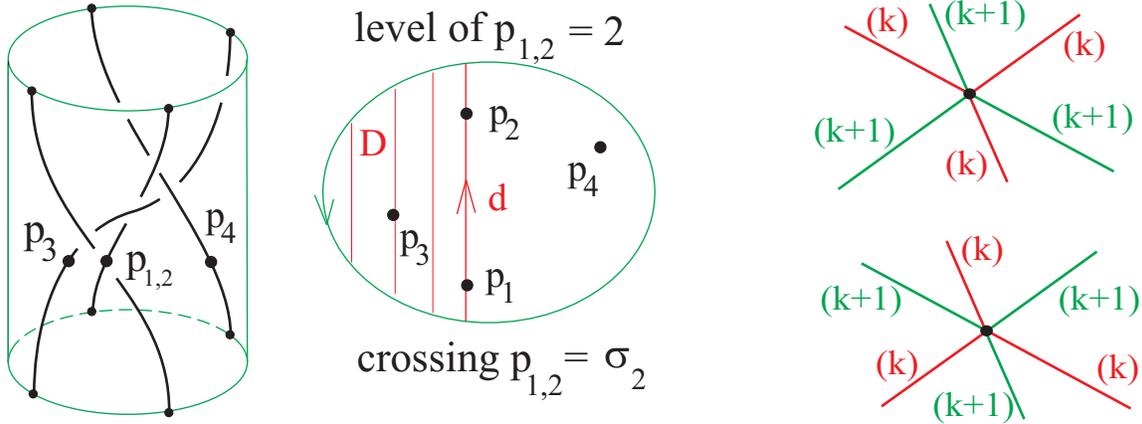}
\caption{Levels in trace graphs of closed braids.}
\label{fig:LevelTraceSubgraphs}
\end{figure}

\begin{definition}
\label{def:LevelsInTraceGraph}
Any point $p\in\tg(\hat\be)$ that is not a vertex
 corresponds to an ordered pair
 $(p_i,p_j)\subset\hat\be\cap(\dxy\times\{z\})$.
Let $t$ be the time moment when $p_i,p_j$ project
 to the same point under $\pr_{xz}:\hat\be\to\axz$.
Then $\rot_t(\hat\be)$ in a thin slice
 $\dxy\times(z-\e,z+\e)$ looks like a braid generator $\si_k$
 or $\si_k^{-1}$, all other strands do not cross each other,
 see Figure~\ref{fig:LevelTraceSubgraphs}.
The index $k$ is called \emph{the level} of the point $p\in\tg(\hat\be)$.
\ed
\end{definition}

>From another point of view we may compute the level
 of a crossing $(p_i,p_j)$ as follows.
Take the oriented straight segment $d$ having endpoints on
 $\bd\dxy\times\{z\}$ and passing through
 $p_i$ first and $p_j$ after.
Complete $d$ with the arc of $\bd\dxy\times\{z\}$
 to get an oriented circuit bounding a disk $D$, 
 see Figure~\ref{fig:LevelTraceSubgraphs}.
The number of intersections $\hat\be\cap\inp D$ plus 1
 is called the \emph{level} $(k)$ of the point $p$.
We chose the name \emph{level} since any crossing
 of $\pr_{xz}(\hat\be)$ is located at its horizontal level
 with respect to $\xa$.
\smallskip

\begin{lemma}
\label{lem:ChangingLevels}
Going along a trace circle of a trace graph, the level of a point $p$ 
 may change only at a triple vertex as follows: $k\mapsto k\pm 1$, 
 see Figure~\ref{fig:LevelTraceSubgraphs}.
\end{lemma}
\begin{proof}
The number of intersections of $\hat\be$
 with the disk $D$ from Definition~\ref{def:LevelsInTraceGraph}
 remains invariant until the segment $d$ passes
 through other points of
 $\hat\be\cap(\dxy\times\{z\})$ apart from $p_i,p_j$ defining
 $p\in\tg(\hat\be)$.
While $p$ passes through a triple vertex of $\tg(\hat\be)$,
 the segment $d$ intersects exactly one strand of $\be$,
 hence the number of points $\hat\be\cap D$ changes by $\pm 1$.
This also follows from $\si_k\si_{k+1}\si_k=\si_{k+1}\si_k\si_{k+1}$.
\end{proof}

Orient the 2---dimensional torus $\zt=S_z^1\times S_t^1$
 in such a way that the first direction is \emph{vertical} along $S_z^1$
 and the second one is \emph{horizontal} opposite to $S_t^1$.

\begin{definition}
\label{def:LevelSubgraphs}
Let $G$ be the trace graph of a closed braid $\hat\be$,
 where $\be\in B_n$.
For each $k=1,\dots,n-1$, denote by $S^{(k)}$
 the $k$-th \emph{level subgraph} consisting of all edges
 having the level $(k)$.
Orient each edge in $G$ vertically along $S_z^1$.
A \emph{right attractor} is an oriented cycle
 $\ra^{(k)}\subset S^{(k)}$ such that at each triple vertex,
 where two edges of $S^{(k)}$ go up,
 the cycle $\ra^{(k)}$ goes to the right.
Denote by $(q^{(k)},r^{(k)})$ \emph{the winding numbers} of
 $\ra^{(k)}$ in the vertical direction $S_z^1$ and
 reversed horizontal direction $(-S_t^1)$, respectively.
Let $e^{(k)}:S^{(k)}\to\zt$ be the $k$-th \emph{level embedding}
 induced by the torus projection
 $\pr_{zt}:S^{(k)}\subset G\to\zt=S_z^1\times S_t^1$,
 see Lemma~\ref{lem:PositionLevelSubgraphs}b below.
\ed
\end{definition}

One right attractor of each $S^{(k)}$, $k=1,2,3$, is shown 
 by fat arcs in Figure~\ref{fig:SplittingIntoLevelSubgraphs}.
In both pictures the 6 marked right attractors have the same winding numbers $(1,0)$.

\begin{lemma}
\label{lem:PositionLevelSubgraphs}
Let $G\subset\T$ be the trace graph of
 a closed braid $\hat\be$, where $\be\in B_n$.
\smallskip

\noindent
{\bf (a)}
Any level subgraph $S^{(k)}$ has only trivalent vertices;
 at each vertex one edge goes down, two edges go up
 or vice versa with respect to the projection $\pr_z:G\to S_z^1$.
\smallskip

\noindent
{\bf (b)}
Any level subgraph $S^{(k)}$ projects 1-1
 to its image under $\pr_{zt}:S^{(k)}\to\zt$.
\smallskip

\noindent
{\bf (c)}
Subgraphs $S^{(k)}$ and $S^{(m)}$ have common points
 under $\pr_{zt}$ if and only if $|k-m|=1$;
 the adjacent subgraphs can meet only in triple vertices
 as in Figure~\ref{fig:LevelTraceSubgraphs}.
\smallskip

\noindent
{\bf (d)}
If $k>m+1$ then the edges of $S^{(k)}$ overcross
 $S^{(m)}$ under $\pr_{zt}:G\to\zt$.
\smallskip

\noindent
{\bf (e)}
Each level subgraph $S^{(k)}$ has at least
 one right attractor.
Its vertical winding number $q^{(k)}$ is positive.
Any two right attractors in $S^{(k)}$ have no
 common points.
\smallskip

\noindent
{\bf (f)}
Under the shift $t\mapsto t+\pi$ each level subgraph $S^{(k)}$
 maps to the subgraph $S^{(n-k)}$.
\end{lemma}
\begin{proof}
{\bf (a)}
A triple vertex $v\in G$ corresponds to
 a triple intersection $(qrs)$ of strands from $\be$,
 see Figure~\ref{fig:MarkingsTraceCircles}.
Let $k$ be the level of the crossing $p$ formed
 by the distant strands of $q$ and $s$ right below $(qrs)$.
By Lemma~\ref{lem:PositionTraceCircles}a the crossing $p$ is associated to 
 a point below $v$ in the middle trace circle passing through $v$.
Right above $p$ the crossings formed by
 the strands $(qr)$ and $(rs)$ have the same level $k$.
The three remaining types of crossings have
 the same level $k+1$ or $k-1$, see Figure~\ref{fig:LevelTraceSubgraphs}.
\medskip

\noindent
{\bf (b)}
If the trace graph $G$ has a crossing under
 the projection $\pr_{zt}:G\to\zt$ then
 the points forming the crossing have the same $z$-coordinate
 and different $x$-coordinates.
Hence they correspond to 2 crossings of some
 diagram $\pr_{xz}(\rot_t(\hat\be))$.
Definition~\ref{def:LevelsInTraceGraph} implies that the levels 
 of these crossings differ at least by 2.
\smallskip

The items {\bf (c)} and {\bf (d)} follow directly
 from the above arguments, see Figure~\ref{fig:SplittingIntoLevelSubgraphs}.
\smallskip

\begin{figure}
\includegraphics[scale=1.0]{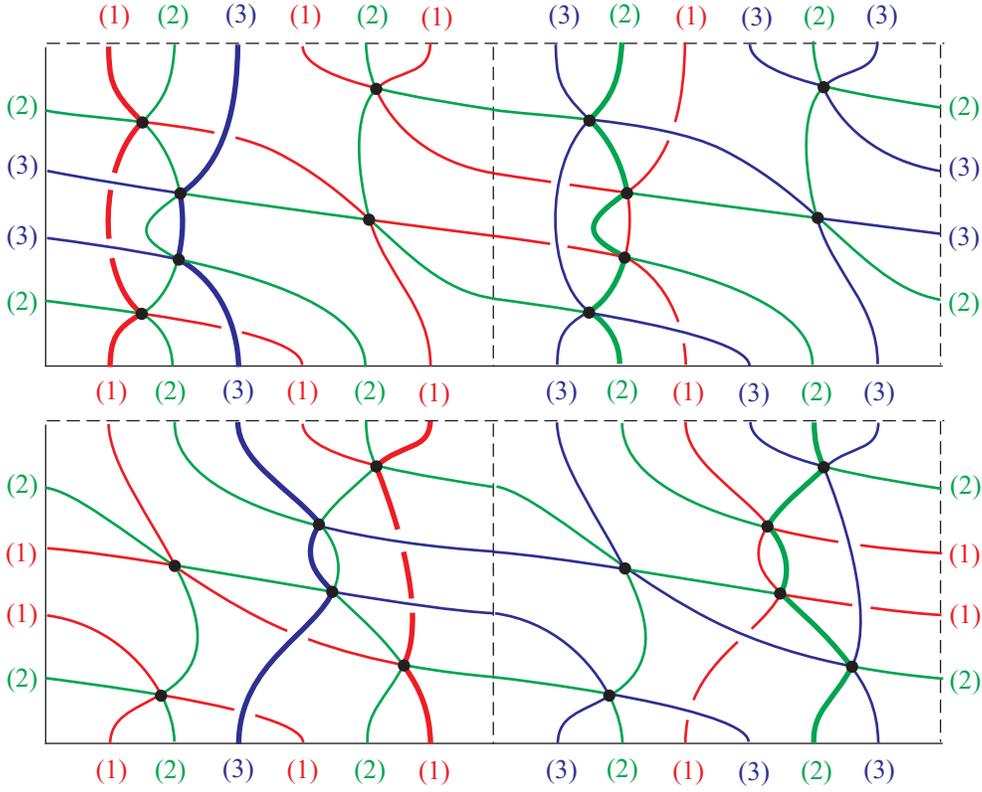}
\caption{Splittings of $\tg(\widehat{\si_2\si_3^2\si_2})$ and
 $\tg(\widehat{\si_2\si_1^2\si_2})$ into level subgraphs.}
\label{fig:SplittingIntoLevelSubgraphs}
\end{figure}

\noindent
{\bf (e)}
Starting with any vertex in $S^{(k)}$ and
 going always to the right in finitely many steps
 we will get a closed cycle oriented vertically, ie $q^{(k)}>0$.
If two right attractors in $S^{(k)}$ have a common vertex then
 they go along the same path and coincide.
\smallskip

\noindent
{\bf (f)}
Let a point $p\in S^{(k)}$ correspond to a pair
 $(p_i,p_j)\in\hat\be$ in a meridional disk
 $\dxy\times\{z\}\subset V$.
The level $k$ is equal to 1 plus the number of
 intersections $\inp D\cap\hat\be$, see Figure~\ref{fig:LevelTraceSubgraphs}.
Under $t\mapsto t+\pi$, the pair $(p_i,p_j)$
 converts to $(p_j,p_i)$, the disk $D$ goes
 to the complementary disk $D'=\dxy\times\{z\}-D$.
Then the level of $(p_j,p_i)$ is 1 plus the number of
 intersections $\inp D'\cap\hat\be$, i.e. $1+(n-2-(k-1))=n-k$.
\end{proof}


\section{Combinatorial encoding trace graphs up to isotopy}
\label{sect:EncodingTraceGraphs}


\subsection{Reconstructing a closed braid using its trace graph}
\label{subs:Reconstruction}
\noindent
\smallskip

The moves on trace graphs are in Figure~\ref{fig:TrihedralMove}
 and Figure~\ref{fig:TetrahedralMoves}.
The trace graph of a closed braid in general position 
 has the combinatorial features summarized below.

\begin{definition}
\label{def:EquivalenceTraceGraphs}
An embedded finite graph $G\subset\T$ is 
 a \emph{generic} trace graph if
\smallskip

$\bu$
under $t\mapsto t+\pi$ the graph $G$
 maps to its image under the symmetry in $S_z^1$;
\smallskip

$\bu$
$G$ splits into \emph{trace circles} monotonic with respect to
 $\pr_z:G\to S_z^1$, they should\\
 \hspace*{8mm}
 intersect in \emph{triple vertices} of $G$ and verify 
 Lemmas~\ref{lem:SplittingTraceCircles}, \ref{lem:PositionTraceCircles};
\smallskip

$\bu$
$G$ splits into $n-1$ \emph{level subgraphs}
 satisfying the conclusions of Lemma~\ref{lem:PositionLevelSubgraphs}.
\medskip

\noindent
A smooth family of trace graphs $\{G_s\}$, $s\in[0,1]$,
 is called \emph{an equivalence} if
\smallskip

$\bu$
 for all but finitely many moments $s\in[0,1]$,
 the trace graphs $G_s$ are generic;
\smallskip

$\bu$
 at each critical moment, $G_s$ changes by a trihedral or
 tetrahedral move.
\medskip

\noindent
An \emph{isotopy} of trace graphs is an equivalence
 through generic trace graphs only.
\ed
\end{definition}

Now we reconstruct a closed braid from its generic trace graph with markings.

\begin{lemma}
\label{lem:Reconstruction}
For a braid $\be\in B_n$,
 the closure $\hat\be\subset V$ can be reconstructed up to isotopy
 in the solid torus from its generic trace graph $\tg(\hat\be)$ with markings.
\end{lemma}
\begin{proof}
Consider a vertical section $P_t=G\cap(\axz\times\{t\})$
 not containing vertices of $G$.
Then $P_t$ is a finite set of points with markings $(ij)[k]$, 
 where $k\in\{1,\dots,\gcd(n_i,n_j)\}$, see Lemma~\ref{lem:SplittingTraceCircles}a.
The points of $P_t$ will play the role of crossings of a diagram of $\hat\be$.
\smallskip

The labelled set $P_t$ defines the \emph{Gauss} diagram $\gd_t$ as follows.
Take $\sqcup_{i=1}^m S_i^1$, split each oriented circle $S_i^1$ into $n_i$ arcs 
 and number them by $1,\dots,n_i$ according to the orientation.
Mark several points in the $q$-th arc of $S_i^1$ in a 1-1 correspondence
 and the same order with the points of $P_t$ projected under $\pr_z:P_t\to S_z^1$
 and having labels $(ij)[k]$ or $(ji)[k]$ 
 for $k=1,\dots,\gcd(n_i,n_j)$.
\smallskip

So each point of $P_t$ gives 2 marked points in $\sqcup_{i=1}^m S_i^1$
 labelled with $(ij)[k]$ and $(ji)[k]$.
Connect them by a chord and get the Gauss diagram $\gd_t$.
The zero Gauss diagram $\gd_0$ is realizable by
 the given diagram of the closed braid $\hat\be$.
The Gauss diagram $\gd_t$ gives rise to a diagram of a closed braid 
 isotopic to $\hat\be$ since the trihedral and tetrahedral moves
 are realized by an isotopy of closed braids, see
 Figure~\ref{fig:MovesOnBraids}.
\end{proof}

Using Lemma~\ref{lem:Reconstruction}, we state
 Theorem~\ref{thm:TraceGraphsMoves} in a slightly different form.

\begin{proposition}
\label{prop:EquivalenceTraceGraphs}
Closed braids $\hat\be_0$ and $\hat\be_1$ are isotopic
 in the solid torus $V$ if and only if
 $\tg(\hat\be_0)$ and $\tg(\hat\be_1)$ are equivalent in
 the sense of Definition~\ref{def:EquivalenceTraceGraphs}.
\qed
\end{proposition}
\smallskip


\subsection{Trace codes of trace graphs}
\label{subs:TraceCodesTraceGraphs}
\noindent
\smallskip

Any curve in $\zt=S_z^1\times S_t^1$ has a homology class $(u,w)$, where
 $u$ is the winding number in the vertical direction $S_z^1$,
 $w$ is the winding number in the direction opposite to $S_t^1$.
Take a generic trace graph $G$ from Definition~\ref{def:EquivalenceTraceGraphs}.

\begin{definition}
\label{def:CyclesHomologyClasses}
A cycle in a level subgraph $S^{(k)}\subset G$ is called \emph{trivial}
 if it bounds a disc under the embedding $e^{(k)}:S^{(k)}\to\zt$.
Any trivial cycle has an orientation induced by
 the oriented torus $\zt$.
Any non-trivial cycle can be oriented in such a way
 that its vertical (possibly, horizontal too)
 winding number is non-negative.
\medskip

\noindent
A level subgraph is said to be \emph{degenerate}, if
 all its non-trivial cycles have homology classes
 that are multiples of each other in $H_1(\zt)=\Z\oplus\Z$.
Let a level subgraph $S^{(k)}$ be non-degenerate.
Denote by $(q^{(k)},r^{(k)})$ the homology class of
 a right attractor.
Among all non-trivial cycles in $S^{(k)}$ choose
 \emph{maximal} cycles with homology classes $(u,w)$ such that
 the value $M=\dfrac{u}{q^{(k)}}-\dfrac{w}{r^{(k)}}$
 is non-zero and maximal.

Recall that $q^{(k)}>0$ by Lemma~\ref{lem:PositionLevelSubgraphs}e.
If $r^{(k)}=0$, then put $M=w$.
The non-degenerate graph $S^{(k)}$ should contain
 non-trivial cycles with $M\neq 0$.
If there are maximal cycles with different homology classes,
 then take one with maximal vertical number $u$.
Now the \emph{maximal} homology class $(u^{(k)},w^{(k)})$ of $S{(k)}$
 is well-defined.
\ed
\end{definition}
\smallskip

By Lemma~\ref{lem:SplittingTraceCircles} trace circles in a trace graph
 are distinguished by their markings.
Any right attractor can be oriented in such a way that
 its vertical winding number is positive.
So right attractors are encoded by cyclic words of vertices.
\smallskip

\begin{definition}
\label{def:TraceCodes}
Choose a base point in each trace circle of a generic trace graph $G$.
Enumerate all vertices of a trace circle $T_{(ab)[i]}\subset G$
 by $(ab)[i]_1,(ab)[i]_2,\dots$
A triple vertex $v\in G$ can be encoded by an
 \emph{ordered triplet}
 $\{ (ab)[i]_x^{(k)}, (ac)[j]_y^{(k\pm 1)}, (bc)[m]_z^{(k)} \}$.
The \emph{trace code} $\tc$ contains
 the following 3 \emph{pieces} of data.
\medskip

\noindent
$\bu$
The \emph{first piece} consists of the ordered triplets associated to
 the vertices of $G$.
\medskip

\noindent
$\bu$
The \emph{second piece} contains the homology classes
 $(q^{(k)},r^{(k)})$ of right attractors\\
\hspace*{3mm}
 for each level subgraph $S^{(k)}$, $k=1,\dots,n-1$.
\medskip

\noindent
$\bu$
The \emph{third piece} is the set of maximal homology classes
 $(u^{(k)},w^{(k)})$ introduced\\
\hspace*{3mm}
 for each level subgraph $S^{(k)}$
 in Definition~\ref{def:CyclesHomologyClasses}, $k=1,\dots,n-1$.
\medskip

\noindent
Two trace codes are called \emph{identical}: $\tc_1=\tc_2$
 if theirs three pieces coincide.
\ed
\end{definition}
\smallskip

Our aim is to reconstruct the embedding of
 a generic closed trace graph $G$ into the thickened torus $\T$
 from its trace code $\tc(G)$, see Lemma~\ref{lem:ReconstructTraceGraph}.
Lemma~\ref{lem:ReconstructLevelSubgraph} proves this for a level subgraph $S^{(k)}\subset G$.
Recall that an \emph{isotopy} in the torus $\zt$
 is a smooth family of diffeomorphisms $F_s:\zt\to\zt$,
 where $s\in[0,1]$, $F_0=\id_{\zt}$.
\smallskip

\begin{lemma}
\label{lem:ReconstructLevelSubgraph}
Let $G$ be the trace graph of a closed braid $\hat\be$,
 where $\be\in B_n$.
\smallskip

\noindent
{\bf (a)}
The embedding $e^{(k)}:S^{(k)}\to\zt$ of
 a \emph{degenerate} level subgraph $S^{(k)}\subset G$
 can be reconstructed by its ordered triplets and
 the homology class $(q^{(k)},r^{(k)})$ of its right attractor
 up to Dehn twists around a right attractor and isotopy in $\zt$.
\smallskip

\noindent
{\bf (b)}
The embedding $e^{(k)}:S^{(k)}\to\zt$ of
 a \emph{non-degenerate} level subgraph can be reconstructed
 up to isotopy in $\zt$ by its ordered triplets,
 the homology class $(q^{(k)},r^{(k)})$ of its right attractor and
 the maximal homology class $(u^{(k)},w^{(k)})$ of $S^{(k)}$.
\end{lemma}
\begin{proof}
{\bf (a)}
A right attractor $\ra^{(k)}\subset S^{(k)}$ can be recognized
 using the set of ordered triplets of vertices.
Embed $\ra^{(k)}$ into $\zt$ according to its
 winding numbers $(q^{(k)},r^{(k)})$.
Add other vertices and edges of $S^{(k)}$ to get an embedding of
 the connected component of $S^{(k)}$ containing the chosen attractor.
If $S^{(k)}$ is non-connected, there is another right attractor
 with the same homology class $(q^{(k)},r^{(k)})$.
\smallskip

We repeat the above steps for all connected components of $S^{(k)}$.
The image of the resulting embedding is contained in
 one or several annuli with the prescribed
 homology class $(q^{(k)},r^{(k)})$.
The whole embedding $S^{(k)}\to\zt$ is well-defined
 up to Dehn twists around a right attractor and isotopy in $\zt$.
\medskip

\noindent
{\bf (b)}
For a non-degenerate subgraph $S^{(k)}$, we construct
 an embedding $S^{(k)}\subset\zt$ as in {\bf (a)}.
We have to improve this embedding by a suitable Dehn twist
 around a right attractor in such a way that
 the maximal homology class is $(u^{(k)},w^{(k)})$.
\smallskip

The number of different homology classes is linear
 with respect to the number of vertices in $S^{(k)}$.
We look at non-trivial cycles in the constructed embedding.
Let $J$ be the algebraic intersection number of
 a right attractor $\ra^{(k)}\subset S^{(k)}$ and
 a non-trivial cycle with a homology class $(u,w)$.
\smallskip

The Dehn twist around $\ra^{(k)}$ acts on
 the homology: $(u,w)\mapsto(u+Jq^{(k)},w+Jr^{(k)})$.
Then $M=\dfrac{u}{q^{(k)}}-\dfrac{w}{r^{(k)}}$
 is invariant under all Dehn twists around $\ra^{(k)}$.
In the already embedded graph $S^{(k)}\subset\zt$
 we may recognize all non-trivial \emph{maximal} cycles with
 the maximal value $M$ computed using $(u^{(k)},w^{(k)})$.
\smallskip

If there are two maximal cycles with different classes
 $(u,w)$ and $(u',w')$, then
 $\dfrac{u-u'}{q^{(k)}}=\dfrac{w-w'}{r^{(k)}}=i$,
 hence $(u,w)=(u',w')+i(q^{(k)},r^{(k)})$ for some $i$.
Since $q^{(k)}$ and $r^{(k)}$ are coprime then $i$ is integer.
Then both cycles have same intersection number $J$
 with the right attractor $\ra^{(k)}$.
So a Dehn twist around $\ra^{(k)}$ acts on
 the set of the homology classes
 of all maximal cycles as a shift by $J(q^{(k)},r^{(k)})$.
\smallskip

We know that among maximal cycles we can find one
 with the homology class $(u^{(k)},w^{(k)})$,
 the vertical number $u^{(k)}$ is maximal possible.
Let $(u,w)$ be the homology class of a maximal cycle $C$
 with the maximal vertical number $u$ in the embedding
 $S^{(k)}\subset\zt$.
There is an integer $i$ such that
 $(u^{(k)},w^{(k)})-(u,w)=iJ(q^{(k)},r^{(k)})$.
\smallskip

The $i$ Dehn twists around the right attractor $\ra^{(k)}$
 convert the cycle $C$ into a required cycle $\ti C$
 with the maximal class $(u^{(k)},v^{(k)})$.
The final embedding $S^{(k)}\subset\zt$ contains
 a basis consisting of $\ti C$ and $\ra^{(k)}$
 with the prescribed homology classes.
Therefore the embedding is well-defined up to isotopy in $\zt$.
\end{proof}
\smallskip

If all level subgraphs of $G$ are non-degenerate,
 we may forget about levels in the trace code $\tc(G)$.
The subgraphs $S^{(k)}\subset G$ should be connected
 and two adjacent subgraphs meet at each triple vertex,
 see Lemma~\ref{lem:PositionLevelSubgraphs}c.
Hence the levels of subgraphs can be reconstructed
 up to the inversion $(1,2,\dots,n-1)\mapsto(n-1,\dots,2,1)$,
 which corresponds to the time shift $t\mapsto t+\pi$.
In the second and third pieces of $\tc(G)$
 we may leave only the homology class of a right attractor
 $\ra^{(1)}$ and the maximal homology class $(u^{(1)},w^{(1)})$
 of the first level subgraph $S^{(1)}$ only.
\smallskip


\section{Recognizing trace graphs in polynomial time}
\label{sect:RecognizingTraceGraphs}


\subsection{Recognizing trace graphs up to isotopy}
\label{subs:RecognizingUpToIsotopy}

\begin{lemma}
\label{lem:ReconstructTraceGraph}
Two generic trace graphs $G_0$ and $G_1$
 are isotopic in the thickened torus $\T$ if and only if
 their trace codes $\tc(G_0)$ and $\tc(G_1)$ become identical
 after suitable cyclic permutations of vertices in trace circles.
\end{lemma}
\begin{proof}
The part \emph{only if} follows from the fact that
 the trace code is invariant under isotopy in $\T$.
The part \emph{if} says that the embedding of a trace graph $G$
 into the thickened torus $\T=\axz\times S_t^1$
 can be reconstructed from its trace code.
\smallskip

By Lemma~\ref{lem:ReconstructLevelSubgraph} we may reconstruct embeddings of level
 subgraphs $S^{(k)}\subset G$ into the torus $\zt$.
Two embeddings of $S^{(1)}$ and $S^{(2)}$ can be
 joint together since the union $S^{(1)}\cup S^{(2)}$
 should be embedded into $\zt$ by Lemma~\ref{lem:PositionLevelSubgraphs}c.
The resulting embedding is well-defined up to isotopy
 in $\zt$ provided that either one of the subgraphs
 $S^{(1)}$ and $S^{(2)}$ is non-degenerate or
 their right attractors have distinct homology classes.
\smallskip

We embed the third subgraph $S^{(3)}$ into $\zt$ to get
 a joint embedding $S^{(2)}\cup S^{(3)}\subset\zt$ as above.
The union $S^{(1)}\cup S^{(2)}\cup S^{(3)}$ can be already
 considered as an embedding into the thickened torus $\T$
 since the edges of $S^{(3)}$ should overcross $S^{(1)}$ in $\zt$.
\smallskip

The final embedding $G\subset\T$ is well-defined up to isotopy
 if either one of the subgraphs $S^{(k)}$ is non-degenerate or
 there are two right attractors with different homology classes.
Otherwise all $S^{(k)}$ are degenerate and
 the embedding $G\subset\T$ is invariant under
 3-dimensional Dehn twists around the common right attractor.
\end{proof}
\smallskip

Proposition~\ref{prop:RecognizeUpToIsotopy} gives a (surprisingly) polynomial algorithm
 recognizing complicated topological objects: 
 trace graphs up to isotopy in a thickened torus.
\smallskip

\begin{proposition}
\label{prop:RecognizeUpToIsotopy}
Let $\be,\be'\in B_n$ be braids of length $\leq l$.
There is an algorithm of complexity $C(n/2)^{n^2/8}(6l)^{n^2-n+1}$
 to decide whether $\tg(\hat\be)$ and $\tg(\hat\be')$
 are isotopic in the thickened torus $\T$,
 where the constant $C$ does not depend on $l$ and $n$.
In the case of pure braids, the power $n^2/8$ can be replaced by 1.
If the closure of a braid is a \emph{knot}, a single circle in the solid torus,
 then the complexity reduces to $Cn(6l)^{n-1}$.
\end{proposition}
\begin{proof}
By Lemma~\ref{lem:AlgConstruction} we may assume that the trace graphs
 $\tg(\hat\be),\tg(\hat\be')$ have $Q\leq 2l(n-2)$ triple vertices.
If we fix numbers $k$ in markings $(ij)[k]$ with $i\neq j$ and
 a base point in each trace circle then
 we can construct trace codes $\tc(\hat\be),\tc(\hat\be')$
 of $\tg(\hat\be)$, $\tg(\hat\be')$, see Definition~\ref{def:TraceCodes}.
The trace codes $\tc(\be),\tc(\be')$ can be compared
 in linear time with respect to the number $Q$ of triple vertices.
\smallskip

By Lemma~\ref{lem:SplittingTraceCircles}a the graph $\tg(\hat\be)$ 
 splits into $N(\be)$ trace circles.
Denote by $k_1,\dots,k_{N(\be)}$ the number of triple vertices
 in the trace circles of $\tg(\hat\be)$.
Then there are exactly $k_1k_2\cdots k_{N(\be)}$ choices of
 base points in the trace cirles.
Since $k_1+\cdots+k_{N(\be)}=3Q\leq 6l(n-2)$, we have
$$k_1k_2\cdots k_{N(\be)}
 \leq \left(\dfrac{ k_1+\cdots+k_{N(\be)} }{N(\be)} \right)^{N(\be)} 
 \leq\left(\dfrac{6l(n-2)}{N(\be)}\right)^{N(\be)}
 \leq (6l)^{n^2-n}$$ due to the estimates
 $n-1\leq N(\be)\leq n(n-1)$ from Lemma~\ref{lem:SplittingTraceCircles}b.
\smallskip

Let $(n_1,\dots,n_m)$ be the lengths of cycles
 of the induced permutation $\ti\be\in S_n$.
There are $\leq n/2$ \emph{non-trivial} cycles with lengths $n_i>1$.
For each pair of non-trivial cycles with lengths $(n_i,n_j)$, there are
 $\gcd(n_i,n_j)$ choices of numbers $k$ in markings $(ij)[k]$,
 ie totally $\prod_{i<j}\gcd(n_i,n_j)$.
Since the number of pairs is $\leq\binom{n/2}{2}$ and
 $\gcd(n_i,n_j)\leq n/2$, the number of choices
 $\leq(n/2)^{\binom{n/2}{2}}\leq (n/2)^{n^2/8-1}$ for $n\geq 4$.
\smallskip

With a fixed choice of markings and base points, we check whether
 $\tc(\hat\be)=\tc(\hat\be')$ with complexity $Cl(n-2)$.
So, the final complexity of the algorithm
 is $C(n/2)^{n^2/8}(6l)^{n^2-n+1}$.
For pure braids, markings $(ij)$ without $[k]$
 are well-defined and we may replace $n^2/8$ by 1.
If $\hat\be$ is a knot, then $N(\be)=n-1$, markings
 $[k]\in\{1,\dots,n-1\}$ are well-defined and
 the complexity reduces to $Cn(6l)^{n-1}$.
\end{proof}


\subsection{Recognizing trace graphs up to trihedral moves}
\label{subs:RecognizingUpToTrihedralMoves}
\noindent
\smallskip

Now we extend Proposition~\ref{prop:RecognizeUpToIsotopy} 
 to recognize trace graphs up to trihedral moves.

\begin{definition}
\label{def:ReducedTraceGraph}
Let $G$ be a generic trace graph from Definition~\ref{def:EquivalenceTraceGraphs}.
A \emph{trihedron} $T\subset G$ is a subgraph homeomorphic to the graph 
 $\theta$ with 3 edges connecting 2 vertices.
A trihedron in a generic trace graph $G$ is called \emph{embedded} 
 if the interiors of its edges do not contain vertices of $G$.
After eliminating (in any order) all embedded trihedra of
 $\tg(\hat\be)$ we get a \emph{reduced} trace graph
 $\ov{\tg}(\hat\be)$.
\ed
\end{definition}
\smallskip

\begin{lemma}
\label{lem:ReducedTraceGraphs}
Let $\ov{\tg}(\hat\be),\ov{\tg}(\hat\be')$ be
 reduced trace graphs of closed braids $\hat\be,\hat\be'$,
 respectively.
The original trace graphs $\tg(\hat\be),\tg(\hat\be')$
 are equivalent through trihedral moves if and only if
 the reduced graphs $\ov{\tg}(\hat\be),\ov{\tg}(\hat\be')$
 are isotopic in $\T$.
\end{lemma}
\begin{proof}
The part \emph{if} is trivial since reduced graphs
 are obtained by trihedral moves.
\smallskip

The part \emph{only if}.
The given equivalence between the original graphs
 provides an equivalence $\{G_s\}$ through trihedral moves only,
 where $s\in[0,1]$, $G_0=\ov{\tg}(\hat\be)$ and $G_1=\ov{\tg}(\hat\be')$.
The trihedral moves in $\{G_s\}$ can create or delete
 only embedded trihedra.
We simulate the creation of each trihedron $T$ as shown in 
 Figure~\ref{fig:SimulationTrihedron}.

\begin{figure}[!h]
\includegraphics[scale=1.0]{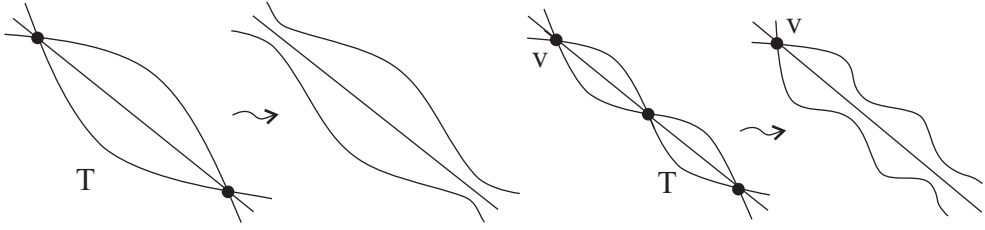}
\caption{Simulation of the appearance of a trihedron.}
\label{fig:SimulationTrihedron}
\end{figure}

Either $T$ will dissappear completely
 by a further trihedral move in $\{G_s\}$ or
 an adjacent trihedron will be deleted and will destroy $T$.
In both cases we miss the deleting move in the simulation.
After simulating all trihedral moves the equivalence
 $\{G_s\}$ becomes a required isotopy between reduced graphs.
\end{proof}
\medskip

\noindent
{\bf Proof of Theorem~\ref{thm:RecognizingUpToTrihedralMoves}.}
Embedded trihedra in a trace graph can be recognized
 in quadratic time with respect to the number of vertices.
For each pair of vertices, we check if they are connected by
 three edges not containing other vertices.
After that the algorithm of Proposition~\ref{prop:RecognizeUpToIsotopy} 
 can be applied to the reduced trace graphs $\ov{\tg}(\hat\be),\ov{\tg}(\hat\be')$
 and gives the required polynomial complexity in the braid length.
\qed
\medskip

A \emph{meridional quadrisecant} of a closed braid $\hat\be$
 in the solid torus $V$ is a straight line in a meridional disk $\dxy\times\{z\}$    
 meeting $\hat\be$ in 4 points.
For an equivalence $\{\hat\be_s\}$ without meridional quadrisecants,
 the canonical loops of rotated braids $\rot_t(\hat\be_s)$ can pass
 only through $\Si_{\tangint},\Si_{\cubicsn}$ and can touch $\Si_{\trip}$, 
 see subsection~\ref{subs:SingularitiesBraids}.
Passing through $\Si_{\quadrup}$ creates a
 a meridional quadrisecant in a closed braid.
Passing through a tangency with $\Si_{\trip}$
 corresponds to a trihedral move in Figure~\ref{fig:TrihedralMove}.
\smallskip

\begin{corollary}
\label{cor:NoQuadrisecants}
Let $\be,\be'\in B_n$ be braids of length $\leq l$.
There is an algorithm of complexity $C(n/2)^{n^2/8}(6l)^{n^2-n+1}$
 to decide whether there is an equivalence $\{\hat\be_s\}$
 such that $\cl(\hat\be_s)$ can pass only through $\Si_{\tangint},\Si_{\cubicsn}$
 and can touch $\Si_{\trip}$ for $s\in[0,1]$.
\end{corollary}
\begin{proof}
The closed braids $\hat\be,\hat\be'$ are equivalent in the above sense
 if and only if their trace graphs $\tg(\hat\be),\tg(\hat\be')$
 are equivalent through trihedral moves only.
So the algorithm of Theorem~\ref{thm:RecognizingUpToTrihedralMoves} can be applied
 to $\tg(\hat\be),\tg(\hat\be')$.
\end{proof}
\smallskip


\section{A geometric recognizing 3-braids up to conjugacy}
\label{sect:Recognizing3Braids}

According to Gonz\'alez-Meneses \cite{GM},
 if two braids $\al$ and $\be$ satisfy $\al^k=\be^k$ in $B_n$
 for some $k\neq 0$, then $\al$ and $\be$ are conjugate.
It follows that braids $\al$ and $\be$ are conjugate
 if and only if $\al^k$ and $\be^k$ are conjugate
 for some $k\neq 0$, see Gonz\'alez-Meneses \cite[Corollary~1.2]{GM}.
For any braid $\be\in B_n$, there is a power $k$ such that 
 the permutation $\ti\be^k\in S_n$ induced by $\be^k$ is trivial, hence $\be^k$ is pure.
So the conjugacy problem for the braid group $B_n$ reduces to the case of pure braids.


\subsection{Cyclic invariants based on 3-subbraids}
\label{subs:CyclicInvariants}
\noindent
\smallskip

In this subsection we recognise closed pure 3-braids up to isotopy
 in the solid torus by using invariants of their trace graphs
 calculable in a linear time with respect to the braid length.
Then trace circles in the trace graph $\tg(\hat\be)$
 can be denoted simply by $T_{(ij)}$, where $i,j\in\{1,\dots,n\}$.
We shall define cyclic invariants depending on 3-subbraids of $\be$ 
 and distinguishing all pure 3-braids up to conjugacy.
\smallskip

Take a pure braid $\be\in B_n$ and enumerate the components of
 $\hat\be$ by $1,\dots,n$.
Fix three pairwise disjoint indices $i,j,k\in\{1,\dots,n\}$.
We shall define the cyclic invariants $C_{(ij)}$
 depending on the 3-subbraid $\be_{ijk}$ based
 on the strands $i,j,k$.
\smallskip

\begin{figure}
\includegraphics[scale=1.0]{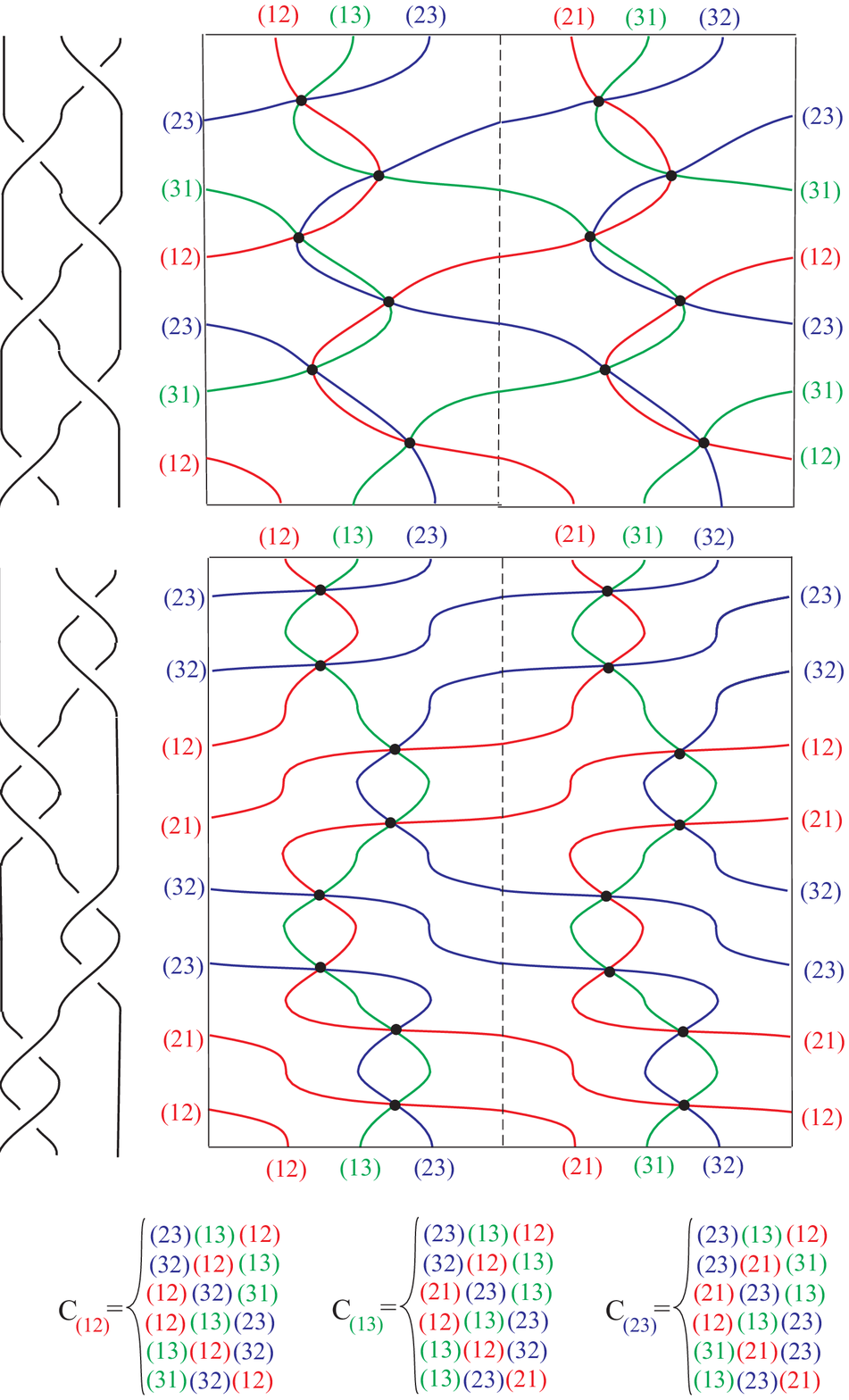}
\caption{Trace graphs of the Borromean links:
 $(\si_1\si_2^{-1})^3$ and $\si_1^2\si_2^2\si_1^{-2}\si_2^{-2}$.}
\label{fig:BorromeanLinks}
\end{figure}

\begin{definition}
\label{def:CyclicInvariants}
Take the reduced trace graph $\ov{\tg}(\hat\be_{ijk})$
 well-defined up to isotopy of $\be_{ijk}$
 by Lemma~\ref{lem:ReducedTraceGraphs} since tetrahedral moves
 are not applicable for 3-braids.
For each triple vertex $v\in T_{(ij)}$, we write
 the ordered triplet of the markings of trace circles passing
 through $v$ in the order from left to right below $v$, see
 Figure~\ref{fig:DynamicTripleVertices}.
The vertices and their triplets are ordered vertically
 in the direction $S_z^1$.
Then $C_{(ij)}(\hat\be)$ is a vertical column of triplets,
 the invariant is defined up to cyclic permutation, 
 see Figure~\ref{fig:BorromeanLinks}.
Similarly we define
 $C_{(ik)}(\hat\be)$, $C_{(jk)}(\hat\be)$.
\ed
\end{definition}
\smallskip

Due to the symmetry of $\tg(\hat\be)$ under the shift $t\mapsto t+\pi$,
 the other invariants $C_{(ji)},C_{(ki)},C_{(kj)}$ can be reconstructed
 from the already defined ones.
\smallskip

\begin{example}
Figure~\ref{fig:BorromeanLinks} contains the trace graphs of the closures of
 the 3-braids $(\si_1\si_2^{-1})^3$ and
 $\si_1^2\si_2^2\si_1^{-2}\si_2^{-2}$.
Both closures are Borromean links, ie the braids are conjugate.
In fact the second graph can be isotoped to the first one
 by eliminating the couple of embedded trihedra.
The cyclic invariants $C_{(12)},C_{(13)},C_{(23)}$
 are shown below the pictures.
The vertices of the embedded trihedra in the second trace graph
 are encoded by $(12)(13)(23)$ and $(23)(13)(12)$,
 ie the extreme markings swap their positions.
Moreover, the cyclic invariants show that the braids are not trivial.
\end{example}
\smallskip


\subsection{Recognizing 3-braids up to conjugacy in a linear time}
\label{subs:ConjugacyInLinearTime}

\begin{lemma}
\label{lem:ComplexityCyclicInvariants}
Number components of two closed pure 3-braids $\be,\be'$ by 1,2,3.
Suppose that the ordered links $\hat\be,\hat\be'$
 are isotopic in the solid torus $V$.
Then the cyclic invariants $C_{(ij)}(\hat\be)$ and
 $C_{(ij)}(\hat\be')$ coincide for all disjoint $i,j\in\{1,2,3\}$.
The invariant $C_{(ij)}(\hat\be)$ is calculable in linear time
 with respect to the length of $\be$.
\end{lemma}
\begin{proof}
By Proposition~\ref{prop:EquivalenceTraceGraphs} the trace graphs of isotopic
 closed braids are connected by an isotopy in the thickened torus $\T$,
 trihedral moves and tetrahedral moves.
The cyclic invariants do not changed under isotopy of trace graphs.
Tetrahedral moves are not applicable for 3-braids.
Trihedral moves create trihedra that are recognizable by
 cyclic invariants and deleted in the construction of
 Definition~\ref{def:CyclicInvariants}.
To compute $C_{(ij)}(\hat\be)$ we need to look at
 all triple vertices of the trace circle $T_{(ij)}$.
The total number of vertices is not more than $2l$ by Lemma~\ref{lem:AlgConstruction}.
\end{proof}
\smallskip

Recall that closed pure 2-braids are classified up to conjugacy
 by the linking number $\lk_{12}$ of closed strands 1 and 2.
Proposition~\ref{prop:CyclicInvariantsComplete} implies that 3-braids can be recognized 
 up to conjugacy in linear time with respect to their length.
\smallskip

\begin{proposition}
\label{prop:CyclicInvariantsComplete}
Fix closed pure 3-braids $\be,\be'$ with ordered components.
The braids $\be,\be'$ are conjugate if and only if
 the linking numbers $\lk_{12}(\hat\be)=\lk_{12}(\hat\be')$ and
 the cyclic invariants $C_{(12)}(\hat\be),C_{(12)}(\hat\be')$
 coincide up to cyclic permutation.
\end{proposition}
\begin{proof}
The part \emph{only if} is Lemma~\ref{lem:ComplexityCyclicInvariants}.
The part \emph{if} says the original 3-braid can be
 reconstructed from its invariants $\lk_{12}$ and $C_{(12)}$.
Simply assume that $\lk_{12}=0$, eg strands 1, 2 are straight,
 ie multiply both braids by $\De^{-2\lk_{12}}$, where
 $\De=(\si_1\si_2)^3$.
\smallskip

Consider a meridional disk $D_z=\dxy\times\{z\}$ in the solid torus $V$,
 where the closed braid $\hat\be$ lives.
Mark the intersection points $D_z\cap\hat\be$ by 1, 2, 3
 according to the components of $\hat\be$.
Since points 1, 2 do not move in $D_z$ while $z$ varies,
 we need to know only how point 3 moves through the line connecting points 1, 2.
There are 6 cases that are 
 distiguished by ordered triplets from $C_{(12)}$.

\begin{figure}[!h]
\includegraphics[scale=1.0]{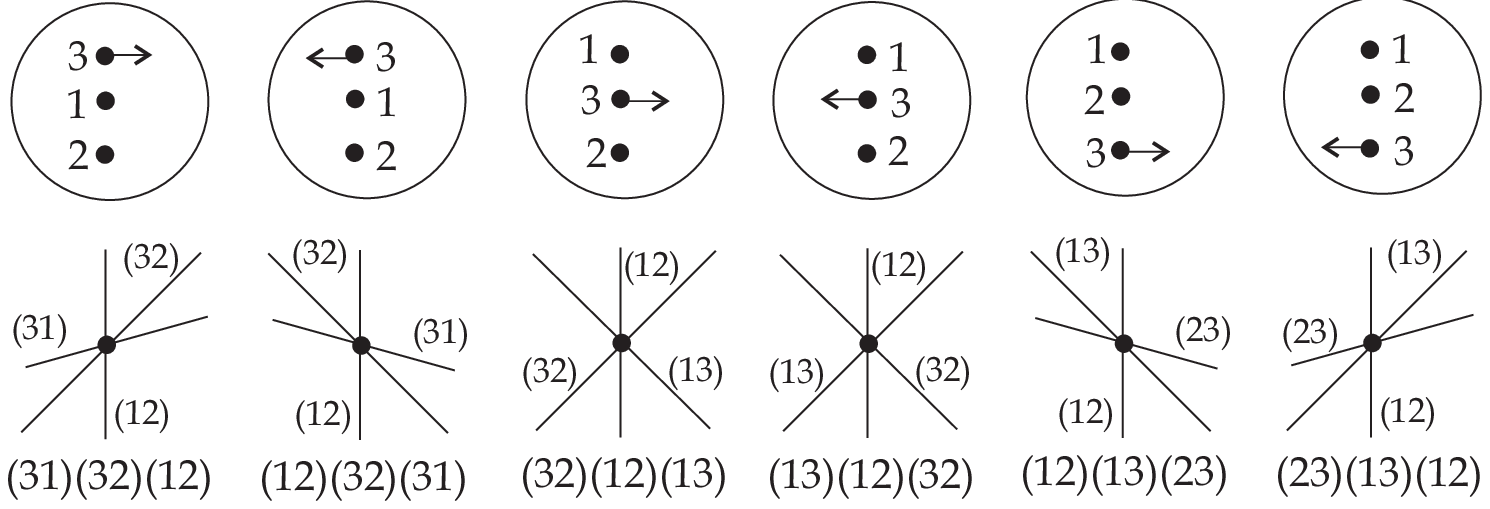}
\caption{Dynamic interpretation of triple vertices.}
\label{fig:DynamicTripleVertices}
\end{figure}

In Figure~\ref{fig:DynamicTripleVertices} the arrow at point 3 
 shows its meridional velocity while $z$ increases.
In the left picture, the line connecting points 1, 3 is going to have
 a positive slope in $D_z$, hence the trace circle $T_{(31)}$
 is increasing as a function $z(t)$.
So triplets describes neighbourhoods of associated triple vertices,
 which can be joined together to get a complete trace graph leading
 to a braid by Lemma~\ref{lem:Reconstruction}.
\end{proof}
\smallskip


\end{document}